\documentclass[12pt,reqno]{amsart}

\usepackage{latexsym}
\usepackage{graphicx}

\newtheorem{theo}{{\sc Theorem}}[section]

\newtheorem{cor}[theo]{{\sc Corollary}}

\newtheorem{lem}[theo]{{\sc Lemma}}
\newtheorem{prop}[theo]{{\sc Proposition}}
\newtheorem{claim}[theo]{{\sc Claim}}

\theoremstyle{remark}
\newtheorem{remark}[theo]{{\sc Remark}}

\def\rr{\mathbb R}
\def\mm{\mathbb M}
\def\nn{\mathbb N}
\def\cc{\mathbb C}

\def\erw{\mathbb E}
\def\hh{\mathbb H}

\def\phi{\varphi }

\def\calc{{\mathcal C}}

\def\cali{{\mathcal I}}
\def\calj{{\mathcal J}}

\def\calm{{\mathcal M}}
\def\caln{{\mathcal N}}

\def\cals{{\mathcal S}}

\def\fra{{\mathfrak a}}

\def\frg{{\mathfrak g}}

\def\frk{{\mathfrak k}}
\def\frl{{\mathfrak l}}

\def\fro{{\mathfrak o}}
\def\frp{{\mathfrak p}}

\def\frs{{\mathfrak s}}
\def\frt{{\mathfrak t}}

\def\frv{{\mathfrak v}}

\newcommand{\R}{\mathbb R}
\newcommand{\N}{\mathbb N}

\def\tra{^{\prime}}

\begin{document}

\title[Large deviations for random matrix ensembles]
{Large deviations for random matrix ensembles in mesoscopic
physics}

\author{Peter Eichelsbacher}
\address{Ruhr-Universit\"at Bochum, Fakult\"at f\"ur Mathematik,
 NA 3/68, D-44780 Bochum, Germany}
\email{peter.eichelsbacher@ruhr-uni-bochum.de}

\author{Michael Stolz}
\address{Ruhr-Universit\"at Bochum, Fakult\"at f\"ur Mathematik, NA 4/32, D-44780 Bochum, Germany}
\email{michael.stolz@ruhr-uni-bochum.de}

\thanks{Both authors have been supported by
Deutsche Forschungsgemeinschaft via SFB/TR 12. They would like to
thank Thomas Kriecherbauer, Margit R\"osler, and Jens Sommerauer for helpful
discussions.}
\date{\today}

\begin{abstract}
In his seminal 1962 paper on the ``threefold way'', Freeman Dyson
classified the spaces of matrices that support the random matrix
ensembles deemed relevant from the point of view of classical
quantum mechanics. Recently, Heinzner, Huckleberry and Zirnbauer
have obtained a similar classification based on less restrictive
assumptions, thus taking care of the needs of modern mesoscopic
physics. Their list is in one-to-one correspondence with the
infinite families of Riemannian symmetric spaces as classified by
Cartan. The present paper develops the corresponding random matrix
theories, with a special emphasis on large deviation principles.
\end{abstract}

\maketitle

Half a century ago, when physicists started to explore the
usefulness of random matrix ensembles for the study of statistical
properties of the spectra of heavy nuclei, their approach was
firmly rooted in the classical framework of quantum mechanics. The
Hamiltonian of a system was replaced by a random hermitian matrix
each realization of which was supposed to commute with the
appropriate unitary symmetries and with the correct ``time
reversals'', i.e.\ certain antiunitary operators. Concretely, the most
general hermitian matrices that commute with time reversal in the
literal sense are real symmetric matrices, whereas another type of
``time reversal'' symmetry leads to quaternion real matrices (see
\cite[Chapter 1]{Forr} or \cite{Haake} for details). In structural
terms, the spaces of hermitian, real symmetric and quaternion real
matrices can be viewed as tangent spaces to or infinitesimal
versions of Riemannian symmetric spaces (RSS) of type A, AI, AII,
respectively. In his landmark article ``The threefold way''
(\cite{Dyson}) of 1962, Dyson proved that any hermitian matrix
which commutes with a group of unitary and ``time reversal''
symmetries reduces to a block matrix, whose blocks are of the
three types described above.

In the 1990s the theoretical condensed matter physicists Altland
and Zirnbauer argued that random matrix models for so-called
mesoscopic normal-superconducting hybrid structures (\cite{AZPRL},
\cite{AZ}, \cite{PADC}) must be taken from the infinitesimal
versions of the symmetric spaces of class B/D, DIII, C, CI (see
the list of classical symmetric spaces in Section \ref{Nambu}
below). What they had in mind are physical systems which typically
consist of a normalconducting quantum dot, i.e., a small metallic
device of extension $1 \mu m$ or less, which is in contact, via
potential barriers, with two superconducting regions. The metallic
quantum dot may or may not be disordered. In the latter case it is
assumed that the geometric shape of the dot is such that the
classical motion of a single electron inside it is chaotic. There
may be a small magnetic flux present and some impurity atoms may
cause spin-orbit scattering. The temperature is so low that the
size of the quantum dot is much smaller than the phase coherence
length of an electron. The latter feature is truly microscopic. On
the other hand, tiny as it is, the quantum dot consists of a huge
number of atoms, which is a distinctively macroscopic feature.
Systems of this kind, displaying both microscopic and macroscopic
features, are called mesoscopic.

We are not aware of an intuitive argument to explain why one
should expect precisely these new ``symmetry classes'' (series of
classical symmetric spaces) to arise in this context. But there is
a more basic question we can try to give some hints about: Given
that Dyson's theorem is a classification result seemingly in the
most general setting quantum mechanics has to offer, why should
there be any leeway for new symmetry classes in mesoscopic systems
as above? The reason is that the standard account of the dynamics
of this kind of system is phrased in the language of second
quantization, i.e., linear superpositions of particle creation and
annihilation operators are acted upon by Hamiltonians. For the
mesoscopic systems described above, however, it is possible to
convert this set-up into something akin to Dyson's first
quantization framework (a Hamiltonian acting on a complex Hilbert
space of ``quasiparticle wavefunctions''). But there remains a
crucial difference: The creation and annihilation operators for
fermions obey the canonical anticommutation relations
$c^{\dagger}_{\alpha} c_{\beta} + c_{\beta} c^{\dagger}_{\alpha} =
\delta_{\alpha \beta},\ c^{\dagger}_{\alpha} c_{\beta}^{\dagger} +
c_{\beta}^{\dagger} c^{\dagger}_{\alpha} = 0,\ c_{\alpha}
c_{\beta} + c_{\beta} c_{\alpha} = 0$, with $\alpha, \beta$ from
an orthonormal basis of single-particle states, and these
relations are mirrored in the translation as a symmetric bilinear
form on the complex Hilbert space (along with the scalar product,
of course). So one has from the outset more structural data than
Dyson had, and it is no surprise that new symmetry classes arise.
(A formalized version of this argument is in the first paragraph
of Section \ref{Nambu} below. For details see the introduction to
\cite{HHZ}).

Let us mention briefly that another source for new symmetry
classes are Dirac fermions in a random gauge field background,
leading to classes AIII, BDI, and CII of the list in Section
\ref{Nambu} below (\cite{Verb}, \cite{Forr}, \cite[Sec.\ 2.3 and
6.2]{HHZ}).

After the physics examples for the new symmetry classes had been
recognized, Zirnbauer and complex geometers Heinzner and
Huckleberry took up the task of updating Dyson's classification
result to the new, enriched framework. In \cite{HHZ}, they proved
that the constituents which make up the Hamiltonian are in
one-to-one correspondence with the ten infinite families of
classical symmetric spaces as classified by Cartan (which
motivates our choice of labels for the
symmetry classes).\\

Like Dyson's ``threefold way'', the ``tenfold way'' of \cite{HHZ}
is established in geometrical terms, without reference to
probability measures on the matrix spaces in question. It is the
object of the present paper to provide a mathematical treatment of
the corresponding random matrix theories. Besides being based on
the systematic framework of \cite{HHZ}, the present article
differs from the existing literature on random matrix ensembles
associated to symmetric spaces (see \cite{Duenez}, \cite{Caselle})
in its focus on Large Deviations Principles. Thus, before going to
business, let us review those aspects of the classical theory  of
the Wigner-Dyson ensembles which will be subsumed  in the present
analysis as instances of symmetry classes A, AI, AII. For
convenience, we only mention the results in the GOE (AI) case. If
$X = (X_{ij})_{1 \le i, j \le n}$ is a symmetric matrix of
real-valued centred Gaussian random variables such that
\begin{itemize}
\item $(X_{ij})_{1 \le i \le j \le n}$ are independent, \item
$\erw(X_{ij}^2) = \frac{1}{2n}\ (i \neq j),\ \erw(X_{ii}^2) =
\frac{1}{n}$,
\end{itemize}
then its distribution is a probability measure on the space of
symmetric $n \times n$ matrices, which is invariant under
conjugation by matrices from the orthogonal group ${\rm O}_n$. The
eigenvalues of $X$, $\lambda_1, \ldots, \lambda_n$, say, are real
valued random variables, and by orthogonal invariance, their joint
distribution has a Lebesgue density $q_n$ that can be given
explicitly:
\begin{equation}
\label{qv1} q_n(x_1, \ldots, x_n)  = \frac{1}{Z_n} \prod_{1 \leq i
< j \leq n}  |x_i - x_j| \exp \bigl( - \frac 12 n \sum_{i=1}^n
x_i^2 \bigr)
\end{equation}
where $Z_n$ is for normalization. As $n \to \infty$, the random
measure
$$ L_n := \frac 1n \sum_{ i = 1}^n \delta_{\lambda_i}$$ tends to a nonrandom limit, namely, to
Wigner's semicircle distribution with density
\begin{equation}
\label{qv50} \frac{1}{\pi} 1_{\{ |x| \leq \sqrt{2} \}} \, \sqrt{2
- x^2}.
\end{equation}
In their paper \cite{BAG}, Ben Arous and Guionnet have carried out
a finer analysis of $L_n$, and have proven that it satisfies a
Large Deviation Principle in $\mathcal{M}_1(\rr)$ (the space of
probability measures on the Borel sets of $\rr$, endowed with the
weak topology) with speed $n^2$ and good rate function
$$
I(\mu) = \frac 14 \int x^2 \, \mu(dx) + \frac 12 \int \int
\log|x-y|^{-1} \,  \mu(dx) \, \mu(dy) - \frac 38,
$$
whose unique minimizer is the semicircle distribution. Recall that
a family of probability measures $(\mu_{\varepsilon})_{\varepsilon
>0}$ on some topological space $X$ is said to obey a Large
Deviation Principle (LDP) with speed $\varepsilon^{-1}$ and good
rate function $I: X\to [0, \infty]$ if
\begin{itemize}
\item $I$ is lower semi-continuous and has compact level sets
$N_L:= \{ x\in X: I(x) \le L\}$, for every $L \in  [0, \infty[$,
\item $\liminf_{\varepsilon \to 0} \varepsilon \log
\mu_{\varepsilon}(G)\ge -\inf_{x\in G} I(x)\quad \forall G
\subseteq X$ open,
 \item $
 \limsup_{\varepsilon \to 0} \varepsilon \log
\mu_{\varepsilon}(A)\le -\inf_{x\in A} I(x)\quad \forall A
\subseteq X$ closed.
\end{itemize}

This paper is organized as follows: In Section \ref{Nambu} we
review the symmetry classification of matrix Hamiltonians in
mesoscopic physics, as given in \cite{HHZ}, and describe the
infinitesimal versions of classical symmetric spaces that turn out
to be in one-to-one correspondence with these symmetry classes.
Then, in Section \ref{qv2}, we introduce probability measures on
these spaces which enjoy invariance properties that reflect those
of quantum mechanical observables. The resulting random matrix
ensembles are called Hamiltonian ensembles. They can be viewed as
generalizations of the Wigner-Dyson ensembles GOE, GUE, GSE. We
use the geometric description of the underlying spaces to derive
the induced joint eigenvalue densities in a uniform way, thus
generalizing (\ref{qv1}) above. In Section \ref{qv46}, we turn to
the large deviations analysis of the empirical eigenvalue measure.
We prove a generalization of the main result of \cite{BAG}, which
covers not only the Hamiltonian ensembles, but
also some matrix ensembles or particle systems of different origin
that have been studied in recent years. In Section \ref{qv43},
then, we describe the Gaussian Hamiltonian ensembles in concrete
terms, analogously to the construction of the GOE which was
reviewed above. We make explicit what the results of Section
\ref{qv46} mean in these special cases.

\section{Symmetries of Nambu space}
\label{Nambu}

Let $(W, b)$ be a complex vector space of dimension $2n\ (n \in
\N)$ together with a nondegenerate symmetric bilinear form $b$. By
polarization we may assume that $W = V \oplus V^*$ and that $b$ is
the natural pairing of $V$ with its dual $V^*$, i.e., $b( x_1 +
\varphi_1, x_2 + \varphi_2) = \varphi_1(x_2) + \varphi_2(x_1)$.
Write $S := \bigwedge V^*$. To $\varphi \in V^*$ assign the wedge
multiplication operator $\epsilon(\varphi) \in \text{\rm End}(S)$,
and to $v \in V$ assign the contraction operator $\iota(v) \in
\text{\rm End}(S)$. The $\epsilon(\varphi)$ and $\iota(v)$ satisfy
the Canonical Anticommutation Relations (CAR) of creation and
annihilation operators on fermionic Fock space, and if $\rho: W
\rightarrow \text{\rm End}(S)$ is given by $\rho(\varphi + v) =
\epsilon(\varphi) + \iota(v)$, then $(\text{\rm End}(S), \rho)$ is
a Clifford algebra $\text{\rm Cliff}(W, b)$ for $(W, b)$. We
regard $W$ as a subspace of the associative algebra $\text{\rm
Cliff}(W, b)$, which we interpret in the usual way as a Lie
algebra. One can embed $\frs\fro(W, b)$ as a Lie subalgebra
consisting of elements which are quadratic in the $w \in W$, and
the adjoint action of $\frs\fro(W, b)$ on $W$ turns out to be
nothing else than the natural action of $\frs\fro(W, b) \subset
\text{\rm End}(W)$ on $W$ (see \cite{GW} for details). In physical
terms, the embedding of $W$ into $\text{\rm Cliff}(W, b) =
\text{\rm End}(S)$ suggests the interpretation of an element of
$W$ as a field operator on fermionic Fock space. The dynamics of a
system of field operators is given by Heisenberg's equation of
motion
$$ i \hbar \frac{dw}{dt} = [H, w],$$
the self-adjoint operator $H$ being the Hamiltonian of the system.
Thus we have seen that if the dynamics of the system is governed
by a quadratic Hamiltonian of a certain type, all relevant
information
is encoded in the structure $(W, b)$.\\

So far we have not yet made explicit that $V$, which plays the
role of the space of single particle states, comes with an
hermitian scalar product $\langle \cdot, \cdot \rangle$. It gives
rise to a ${\mathbb C}$-antilinear bijection $C: V \to V^*: v
\mapsto \langle v,  \cdot \rangle$. $C$ and $\langle \cdot, \cdot
\rangle$ can be extended to the entire space $W$ in such a way
that $\langle w_1, w_2 \rangle = b(C w_1, w_2)\ \forall w_1, w_2
\in W$. The triplet $(W,
b, C)$ is called {\it Nambu space}.\\

Now suppose that a compact group ${\mathcal G}$ acts on $W$ by
unitary or antiunitary transformations. It is inevitable to bring
antiunitary transformations into play here, because time enters
the formalism of quantum mechanics via $i \hbar \frac{d}{dt}$, and
so time reversal is an antiunitary rather than a unitary
transformation. Once $(W, b, C)$, ${\mathcal G}$ and its action on
$W$ are fixed, the translation of Dyson's problem to the present
context is as follows: Describe the space of all hermitian $H \in
\frs\fro(W, b)$ with the property that $Hg = gH$ for all $g \in
{\mathcal G}$. It is convenient to call these $H$ {\it good
Hamiltonians}. Denote by ${\mathcal G}_0$ the group of all
elements of ${\mathcal G}$ which act by unitary transformations.
The basic assumption of \cite{HHZ} is that ${\mathcal G}$ is
generated by ${\mathcal G}_0$ together with at most two elements
acting by antiunitary transformations (One may think of systems
which are invariant under both time reversal and
charge conjugation). \\

It is one of the key insights of Heinzner, Huckleberry and
Zirnbauer in \cite{HHZ}, that it is possible to reduce the problem
to the case  ${\mathcal G}_0 = 1.$ But this comes at the price
that the good Hamiltonians need not be elements of $\frs\fro(W,
b)$, but can belong to $\frs \frp(W, a)$ (for an alternating form
$a$ on $W$) or to $\frs\frl(V)$ (diagonally embedded into
$\text{\rm End}(V) \oplus \text{\rm End}(V^*)$). Write $\frs$ for
any of these three Lie algebras. Although one is ultimately
interested in hermitian operators, one first considers the skew
hermitian elements of $\frs$, which make up a compact real form
$\frg$ of $\frs$. Now let $T$ be an antilinear transformation of
$W$ such that $T^2 = \pm \text{\rm id}$ and suppose that
${\mathcal G} = \langle T \rangle$. From the reduction step in
\cite{HHZ} it emerges that $T$ may be assumed to fix $\frg$. If
$\theta$ denotes conjugation (in $\text{\rm End}(W)$) by $T$, then
$\theta$ restricts to an involutive Lie algebra automorphism of
$\frg$. Let $\frg = \frk \oplus\frp$ be the decomposition into the
$+1$-eigenspace $\frk$ and the $-1$-eigenspace $\frp$ of $\theta$,
the so-called {\it Cartan decomposition}. Then $i\frp$ consists
precisely of those hermitian operators in $\frg$ which commute
with $T$, hence it is the space of good Hamiltonians. If $G$ is a
connected compact Lie group corresponding to $\frg$ and $K$ its
closed subgroup with Lie algebra $\frk$, then $\frp$ can be
thought of as an infinitesimal version of the RSS $G / K$. If
${\mathcal G} = 1$, then the space of good Hamiltonians is $i
\frg$. Since compact Lie groups can be given the structure of an
RSS (see \cite{H1}, Ch. IV \S 6), this case fits into
the overall picture. \\

It is the main result of \cite{HHZ} that the following is the
complete list of spaces of good Hamiltonians that correspond to
Nambu space with the kind of symmetries in question. The labels
refer to Cartan's list of the classical compact Lie
algebras and their involutive automorphisms (see \cite[Ch. X \S
2.3]{H1}), i.e., plainly, to his list of classical
symmetric spaces. Note that in order to obtain the case BDI below
for the full range of parameters, and thus to obtain Cartan's full
list from the symmetries of Nambu space, one has to refine the
above argument in order to take care of the case ${\mathcal G} =
\langle T, T_1\rangle$, $T_1$ being another antilinear
transformation of $W$ with $T_1^2 = \pm \text{\rm id}$. 
We deviate from standard practice in Lie theory in that we do not
require the matrices to be trace-free, in order to recover the
familiar Wigner-Dyson classes as classes A, AI, AII. $X^*$ denotes the conjugate transpose of a complex
matrix $X$.

\begin{description}

\item[Class A] $$  i \frg = \{ X \in {\mathbb C}^{n
\times n}:\ X\ \text{\rm hermitian} \}$$

\item[Class AI] $$  i \frp = \{ X \in {\mathbb R}^{n
\times n}:\ X\ \text{\rm symmetric} \}$$

\item[Class AII] $$  i \frp = \left\{ \left(
\begin{array}{rr} X_1 & X_2 \\ -\overline{X_2} & \overline{X_1}
\end{array} \right):\quad  \begin{array}{l} X_i \in {\mathbb C}^{n \times n}, X_1\ \text{\rm
hermitian},\\ X_2\ \text{\rm skew symmetric} \end{array} \right\}$$

\item[Class AIII] $$  i \frp = \left\{ \left(
\begin{array}{ll} 0 & X \\ X^* & 0 \end{array} \right):\ X \in {\mathbb C}^{s \times t} \right\}$$

\item[Class B/D] $$ i \frg = \{ X \in (i {\mathbb
R})^{n \times n}:\ X\ \text{\rm skew symmetric} \}$$

\item[Class BDI] $$  i \frp = \left\{ \left(
\begin{array}{ll} 0 & X \\ X^* & 0 \end{array} \right):\ X \in (i {\mathbb R})^{s \times t} \right\}$$

\item[Class DIII] $$ i \frp = \left\{ \left(
\begin{array}{rr} X_1 & X_2 \\ X_2 & -X_1
\end{array} \right):\ X_i \in (i {\mathbb R})^{n \times n}\ \text{\rm skew symmetric} \right\}$$

\item[Class C] $$ i \frg = \left\{ \left(
\begin{array}{rr} X_1 & X_2 \\ \overline{X_2} & - \overline{X_1}
\end{array} \right):\quad \begin{array}{l} X_i \in \cc^{n \times n},\ X_1\ {\rm hermitian},\\ X_2\ {\rm symmetric}\end{array}  \right\}$$

\item[Class CI] $$  i \frp = \left\{ \left(
\begin{array}{rr} X_1 & X_2 \\ X_2 & -X_1
\end{array} \right):\ X_i \in {\mathbb R}^{n \times n}\ \text{\rm symmetric} \right\}$$

\item[Class CII] $$  i \frp = \left\{ \left(
\begin{array}{ll} 0 & X \\ X^* & 0 \end{array} \right):\ X \in {\mathbb H}^{s \times t} \right\},$$
where the space $\hh^{s \times t}$ of quaternionic matrices is
embedded into $\cc^{2s \times 2t}$ as
$$ \hh^{s \times t} = \left\{ \left( \begin{array}{rr} U & V \\ - \overline{V} & \overline{U}\end{array}\right):\ U, V \in \cc^{s\times t}\right\}.$$

\end{description}

More precisely, B/D splits into B for $n$ odd, and $D$ for $n$ even.
We will call classes A, AI, AII Wigner-Dyson classes, BDI, AIII and CII chiral classes (in view of their role in modeling Dirac fermions) and the others superconductor or Bogolioubov-de Gennes (BdG) classes (in view of the discussion in the introduction, BdG being a keyword in the conversion to a first 
quantization set-up which was mentioned there).

\section{Hamiltonian ensembles}
\label{qv2}

In this section we randomize the good Hamiltonians, i.e., we put
probability measures on $\tilde{\frg} := i \frg$ resp.\
$\tilde{\frp} := i \frp$. Let $G$ be a connected compact Lie group
with Lie algebra $\frg$, $K$ its closed subgroup with Lie algebra
$\frk$. The adjoint representations ${\rm Ad}_G: G \to {\rm
GL}(\frg)$ and
  ${\rm Ad}_K: K \to {\rm GL}(\frp)$ are given by conjugation of matrices. If $\frv$ is a (nonempty open subset of a)
finite dimensional Euclidian vector space, write ${\rm m}_{\frv}$ for Lebesgue measure on $\frv$. \\
Now we study probability measures on $\tilde{\frg}$ resp.\
$\tilde{\frp}$, restricting our attention to those which are
absolutely continuous w.r.t.\ ${\rm m}_{\tilde{\frg}}$ resp.\
${\rm m}_{\tilde{\frp}}$. Since quantum mechanical observables are
invariant under unitary similarity transformations, it is natural
to assume that the measures are invariant under conjugation with
unitary matrices. So we only consider ${\rm
m}_{\tilde{\frg}}$-densities that are constant on the ${\rm
Ad}_G$-orbits in $\tilde{\frg}$
resp.\ ${\rm m}_{\tilde{\frp}}$-densities that are constant on the ${\rm Ad}_K$-orbits in $\tilde{\frp}$. \\
Now we wish to compute the joint densities of the eigenvalues of a
matrix $X \in \tilde{\frg}$ resp.\ $X \in \tilde{\frp}$ chosen
according to such a measure. These are easy consequences of the
infinitesimal version of Weyl's integration formula for $\frg$ and
its analog, due to Harish-Chandra, for $\frp$. Here we use
standard terminology of elementary Lie theory, see e.g. \cite{OV},
\cite{H1}, \cite{DK}.
To state the
formula for $\frg$, let $T$ be a maximal torus of $G$ with Lie
algebra $\frt$, $W = {\rm N}_G(T)/T $ the corresponding Weyl
group, $R^+ $ and $\frt^+ \subset \frt$ compatible choices of a
system of positive roots and of a positive Weyl chamber, respectively.
\begin{prop}
\label{qv4} There exists $c > 0$ such that for all $f \in {\rm
L}^1(\frg, {\rm m}_{\frg})$ which are constant on $ {\rm
Ad}_G$-orbits there holds
\begin{eqnarray*}
\int_{\frg} f\ d{\rm m}_{\frg} &=& c \int_{\frt^+} f(T) \prod_{\alpha \in R^+} \alpha(T)^2\ {\rm m}_{\frt^+} (dT)\\
&=& \frac{c}{\# W} \int_{\frt} f(T) \prod_{\alpha \in R^+}
|\alpha(T)|^2\ {\rm m}_{\frt} (dT).
\end{eqnarray*}
\end{prop}
\begin{proof}
\cite{DK} Cor. 3.14.2 (ii)
\end{proof}

To state the theorem for $\frp$, write $\fra$ for a maximal
abelian subspace of $\frp$, $W = {\rm N}_K(\fra)/ {\rm C}_K(\fra)$
for the Weyl group of $\frg$ w.r.t.\ $\fra$, $\Sigma^+$ and
$\fra^+ \subset \fra$ for compatible choices of positive
restricted roots (of $\frg$ w.r.t. $\fra$) and of a positive
Weyl chamber, respectively. For $\rho \in \Sigma^+$, $m_{\rho}$ denotes the
multiplicity of $\rho$, i.e. the dimension of the joint
eigenspace, corresponding to the linear form $\rho$, of the
commuting symmetric operators which are induced on $\frg$ by
$\fra$. Alternatively, $m_{\rho}$ is the cardinality of the
inverse image of $\rho$ w.r.t the restriction process described in
\cite{H1}, p. 263-4.

\begin{prop}
\label{qv5} There exists $c > 0$ such that for all $f \in {\rm
L}^1(\frp, {\rm m}_{\frp})$ which are constant on $ {\rm
Ad}_K$-orbits there holds
\begin{eqnarray*}
\int_{\frp} f\ d{\rm m}_{\frp} &=& c \int_{\fra^+} f(A) \prod_{\rho \in \Sigma^+} \rho(A)^{m_{\rho}} {\rm m}_{\fra^+}(dA)\\
&=&  \frac{c}{\# W} \int_{\fra} f(A) \prod_{\rho \in \Sigma^+}
|\rho(A)|^{m_{\rho}} {\rm m}_{\fra}(dA)
\end{eqnarray*}
\end{prop}
\begin{proof}
\cite[Thm. I.5.17]{H2}
\end{proof}
In practical terms, if $f$ is an ${\rm Ad}_G$-invariant ${\rm
m}_{\tilde{\frg}}$-density ${\tilde{\frg}} \to [0, \infty[$ or an
 ${\rm Ad}_K$-invariant ${\rm m}_{\tilde{\frp}}$-density ${\tilde{\frp}} \to [0, \infty[$, then $f(X)$ only depends on $X$ through its eigenvalues.
Now, for all classes except A, AI, AII, the nonzero eigenvalues
come in pairs $\pm \lambda$. Note that in the chiral classes BDI,
AIII and CII, the number of positive eigenvalues (without multiplicity) is $s \wedge t$, 
and later on, when we will let matrix
size tend to infinity, we will have to control the growth of $s =
s(n)$ as well as of $n$. What is more, note that many of the
matrices that make up the spaces of good Hamiltonians are
necessarily of even size, so that there is no question of simply
letting matrix size tend to infinity. Hence we choose the
following framework: $n$ or $(n, s(n))$ is the (pair of) positive integer parameter(s) that will be controlled, $t(n) := n - s(n)$ by
definition. The space of good Hamiltonians for symmetry class
$\calc$ and parameters $n$ or $(n, s(n))$, written as ${\mathbb
M}^{(n)}_{\calc},$ is
contained in $\cc^{d(n) \times d(n)}$. The map
$$ \pi_n:\ \mm_{\calc}^{(n)} \to \overline{\rr^{p(n), +}},$$
where
$$ \rr^{d, +} := \{ x \in \rr^d:\ x_1 > x_2 > \ldots > x_d\},$$
assigns to each $X \in \mm^{(n)}_{\calc} $ the nonincreasing
vector $(\lambda_1, \ldots, \lambda_{p(n)})$ of its (positive) eigenvalues (without multiplicity). By construction, $\pi_n$ separates the adjoint orbits of
$G$ resp.\ $K$, and it is easily seen that $\pi_n(X) \in
\rr^{p(n), +}$ a.s.
\\

Plugging well-known facts about (restricted) root systems (see,
e.g., the appendix to the monograph \cite{OV}) into Propositions
\ref{qv4}, \ref{qv5} above, one obtains the following

\begin{cor}
\label{qv6} Suppose that $X$ is a random element of ${\mathbb
M}^{(n)}_{\calc}$, whose distribution is given by a Lebesgue
density $f$ such that $f = \tilde{f} \circ \pi_n$ for some
measurable $\tilde{f}: \rr^{p(n), +} \to [0, \infty[$.
\begin{itemize}
\item[(a)] If $\calc$ is A, AI or AII, then the joint density of
the eigenvalues of $x$ (in nonincreasing order) is
$$ (x_1, \ldots, x_{p(n)}) \mapsto {\rm const~}\tilde{f}(x_1, \ldots, x_{p(n)}) \prod_{1 \le i < j \le p(n)} (x_i - x_j)^{\beta},$$
where $\beta$ and $p(n)$ are given in the table below. \item[(b)]
Otherwise, the joint density of the positive eigenvalues of $X$
(in nonincreasing order) is given by
$$(x_1, \ldots, x_{p(n)}) \mapsto {\rm const~}\tilde{f}(x_1, \ldots, x_{p(n)}) \prod_{1 \le i < j \le p(n)} (x_i^2 - x_j^2)^{\beta}
\prod_{1 \le i \le p(n)} x_i^{\alpha},$$ where $\alpha, \beta$ and
$p(n)$ are given in the table below.
\end{itemize}
\end{cor}
\vspace{1em}

\begin{center}
\begin{tabular}{|l|l|l|l|l|}
\hline
Class & $d(n)$ & $p(n)$ & $\alpha$ & $\beta$\\
\hline
A & $n$ & $n$ && $2$ \\
AI & $n$ & $n$ & & $1$ \\
AII & $2n$ & $n$ && $4$\\
\hline
BDI & $n$  & $s(n) \wedge t(n)$ & $|s(n) - t(n)|$ & $1$ \\
AIII & $n$ & $s(n) \wedge t(n)$ & $2 |s(n) - t(n)|  + 1$ & $2$\\
CII & $2n$ & $s(n) \wedge t(n)$ & $4 |s(n) - t(n)| + 3$ & $4$\\
\hline
B & $2n + 1$ & $n$ & $2$ & $2$\\
D & $2n$ & $n$ & $0$ & $2$\\
C & $2n$ & $n$ & $2$ & $2$ \\
CI & $2n$ & $n$ & $1$ & $1$\\
DIII ($n$ even) & $2n$ & $\lfloor \frac{n}{2}\rfloor$ & $1$ & $4$\\
DIII ($n$ odd) & $2n$ & $\lfloor \frac{n}{2}\rfloor$ & $5$ & $4$\\
\hline
\end{tabular}
\end{center}
\vspace{2em}

\section{An LDP for the empirical eigenvalue measure}
\label{qv46} In this section we prove a theorem which contains
LDPs for the empirical eigenvalue measures of all Hamiltonian
ensembles, assuming a product structure for $\tilde{f}$ in the
notation of Corollary \ref{qv6}. Note that this subsumes the
well-studied densities on ${\mathbb M}^{(n)}_{\calc}$ of the form
$\exp(- {\rm Tr}(V(X)))$, $V$ a polynomial with positive leading
coefficient (see Remark \ref{qv100} below). This is the set-up:
Let $\beta
>0$, $\gamma \in \N$. Let $\Sigma$ be a closed subinterval of $\R$
if $\gamma$ is odd, of $[0, \infty[$ if $\gamma$ is even. Let
$(w_n)_{n \in \N}$ be a family of continuous nonnegative
real-valued functions on $\Sigma$. For $n \in \nn$ consider random
variables $\Lambda_n = (\lambda_1, \ldots, \lambda_{p(n)})$ with
joint distribution $Q_n = {\mathbb P} \circ \Lambda_{n}^{-1} \in
{\mathcal M}_1(\Sigma^{p(n)})$, given by its Lebesgue density
\begin{equation} \label{qv7}
q_n(x_1, \ldots, x_{p(n)}) = \frac{1}{Z_n} \prod_{1 \leq i < j
\leq p(n)} |x_i^{\gamma} - x_j^{\gamma}|^{\beta}
\prod_{j=1}^{p(n)} w_{n}(x_j)^n \, 1_{\Sigma^{p(n)}} (x_1, \ldots,
x_{p(n)}),
\end{equation}
where $Z_n $ is for normalization.
In what
follows we assume that $p(n) \to \infty$ for $n \to \infty$,
satisfying
\begin{equation} \label{qv8}
\lim_{n \to \infty} \frac{p(n)}{n} = \kappa \in ]0, \infty[.
\end{equation}
Write $\caln(f)$ for the set of zeros of a function $f$. We will
make the following assumptions about the sequence $(w_{n})_n$:
\begin{itemize}
\item[(a1)] there exists a continuous function $w: \Sigma \to [0,
\infty[$ such that
\begin{itemize}
\item $\# \caln(w) < \infty,\ \caln(w_n) \subseteq \caln(w)$ for
large $n$. \item As $n \to \infty$, $w_n \to w$ and $\log w_n \to
\log w$ uniformly on compact sets, \item  $\log w$ is Lipschitz on
compact sets away from $\caln(w)$.
\end{itemize}
\item[(a2)] If $\Sigma$ is unbounded, then there exists $n_0 \in
\nn$ such that
$$ \lim_{x \to \pm \infty} |x|^{\gamma \kappa (\beta \vee 1) + \epsilon} \sup_{n \ge n_0} w_n(x) = 0$$
for some fixed $\epsilon > 0$.
\end{itemize}

For $x = (x_1, \ldots, x_{p(n)}) \in \Sigma^{p(n)}$ set
$$
L_{n}(x) := \frac{1}{p(n)} \sum_{j=1}^{p(n)} \delta_{x_j}.
$$

\begin{theo} \label{qv10}
$\left( {\mathbb P} \circ \left( L_{n} \circ \Lambda_{n}\right)^{-
1}\right)_n = \left( Q_n \circ L_{n}^{-1}\right)_n$ satisfies an
LDP on $\mathcal M_1(\Sigma)$ with respect to the weak topology
with speed $n^2$ and good rate function
\begin{eqnarray} \label{qv11}
I(\mu) & = &\frac{\beta}{2} \kappa^2 \int \int \log |x^{\gamma} - y^{\gamma}|^{-1} \, \mu(dx) \, \mu(dy) \\
& - & \kappa \int \log w(x) \, \mu(dx) - c,
\end{eqnarray}
where $\mu \in \mathcal M_1(\Sigma)$ and
\begin{equation} \label{qv12}
c :=  \lim_{n \to \infty} \frac{1}{n^2} \log Z_n < \infty.
\end{equation}
\end{theo}

\begin{cor} \label{LLN}
If $I$ has a unique minimizer $\mu^*$, then
$$
{\mathbb P} \bigl( L_{n} \circ \Lambda_n  \to \mu^* \bigr) =1
$$
where $\to$ means weak convergence.
\end{cor}

\begin{proof}[Proof of the corollary:]
Using the upper bound of the LDP one obtains the strong law
applying Borel-Cantelli's lemma, see \cite[Theorem II B.3]{Ellis}.
\end{proof}

\begin{remark}
If $\beta \kappa \ge 1$, then it follows from the theory of
logarithmic potentials with external fields, applied to the weight
function $x \mapsto w( x^{1/\gamma})$, that $\mu^*$ exists and
that $I(\mu^*) + c$ is finite, see \cite{ST}, Thm.\ I.1.3 and Ex.\
I.3.5. In this case \eqref{qv12} can be sharpened to $|c| <
\infty$. There exists a vast literature that describes $\mu^*$ in
more detail for various classes of weights, see e.g.\ \cite{ST},
\cite{DKL}, \cite[Thm.\ 3.1]{Erco}. We will restrict ourselves to
giving explicit formulae for $\mu^*$ for Gaussian Hamiltonian
ensembles below in Section \ref{qv43}. In these cases, $\mu^*$ is
``universal" in the sense of Wigner's Theorem. This is the content
of the companion paper \cite{KHC}.
\end{remark}

\begin{remark}
\label{qv100} Although the main focus of the present paper is on
Hamiltonian ensembles associated to (infinitesimal) symmetric
spaces, let us note that Theorem \ref{qv10} contains an LDP for
the much wider class of ensembles which is considered in
\cite[(1.6)]{Erco}. Furthermore, Theorem \ref{qv10} applies to
Wishart matrices (see Remark \ref{qv101} below) and to Jacobi
ensembles, for which the LDP was first proven in \cite{HPjac}. As
observed in \cite{Duenez}, the latter class of ensembles includes
the random matrix ensembles associated to compact symmetric
spaces, except those of classes A, AI, AII.
\end{remark}

The {\it proof} of Theorem \ref{qv10} extends the approach of
\cite{BAG}, \cite{Brasil}, \cite{HPwish}. Define on $\Sigma \times
\Sigma$ the function
$$
F(x,y) := - \frac{\beta}{2} \kappa^2 \log |x^{\gamma} -
y^{\gamma}| - \frac{\kappa}{2} \bigl( \log w(x) + \log w(y)
\bigr),
$$
with $F(x,y) = \infty$ for $x^{\gamma} = y^{\gamma}$ (hence $x =
y$ by definition of $\Sigma$) or $\{ x, y \} \cap {\mathcal N}(w)
\neq \emptyset$, and its truncated versions
$$
F^M(x,y) := F(x,y) \wedge M, \quad M>0.
$$
Moreover, consider the functions
$$
F_n(x,y) := - \frac{\beta}{2} \bigl( \frac{p(n)}{n} \bigr)^2 \log
|x^{\gamma} - y^{\gamma}| - \frac{p(n)}{2n} \bigl( \log w_{n}(x) +
\log w_{n}(y) \bigr)
$$
and their truncated versions
$$
F^M_n(x,y) := F_n(x,y) \wedge M, \quad M>0.
$$
From the definition of $F_n$ it follows that
\begin{equation} \label{Fn}
q_n(x_1, \ldots, x_{p(n)}) = \frac{1}{Z_n} \exp \biggl( - \frac{2
n^2}{p(n)^2} \sum_{1 \leq i < j \leq p(n)} F_n(x_i, x_j) +
\frac{n}{p(n)} \sum_{i=1}^{p(n)} \log w_{n}(x_i) \biggr).
\end{equation}

\begin{lem} \label{qv16}
~
\begin{itemize}
\item[(i)]  For any $M>0$, $F^M_n(x,y)$ converges to $F^M(x,y)$
uniformly as $n \to \infty$.

\item[(ii)] $F$ is bounded from below.
\end{itemize}
\end{lem}

\begin{proof}
The estimate $\log |x-y| \leq \log(|x|+1) + \log(|y|+1)$ implies
\begin{equation}
\label{qv45}
 F_n(x,y) \geq  - \frac{p(n)}{2n} \biggl[ \log \bigl(
(|x^{\gamma}|+1)^{\beta p(n)/n} \, w_n(x) \bigr) +\log \bigl(
(|y^{\gamma}|+1)^{\beta p(n)/n} \, w_n(y) \bigr) \biggr].
\end{equation}
Observe that $\log \bigl( (|x^{\gamma}|+1)^{\beta p(n)/n} \,
w_n(x) \bigr)$ is bounded from above by (a2) and the continuity of
$w_n$. Invoking the full strength of (a2), one sees that for each
$M > 0$ there exist $n_0 \in \nn,\ R_M > 0,\ \delta_{\nu, M} > 0\
(\nu \in \caln(w))$ such that $F_n(x, y) \ge M$ holds for all $n
\ge n_0$ on
$$ A_M := \{ |x| \vee |y| > R_M \} \cup \bigcup_{\nu \in \caln(w) }  \{ |x - \nu| \wedge |y - \nu| <
\delta_{\nu, M} \}.$$ By compactness of $A_M^c$, (i) follows from
(a1) and the definition of $F^M, F_n^M$. One also has $F \ge M$ on
$A_M$. $F$ being continuous, it is bounded on $A_M^c$, and this
proves (ii).
\end{proof}

Our strategy is to first consider the finite positive measures
$P_n := Z_n Q_n$. As to the {\it upper bound}, note that
$$L_{n} \otimes L_{n} (\{(x,y) \in \Sigma^2: x=y\})=
\frac{1}{p(n)}$$ ${\rm m}_{\rr^{p(n)}}$-almost surely, since the
eigenvalues are a.s.\ distinct under the product Lebesgue measure.
Hence almost surely
\begin{equation} \label{vor}
\int \int_{x \not= y} F^M_n(x,y) \, L_{n}(dx) \, L_{n}(dy) = \int
\int F^M_n(x,y) \,  L_{n}(dx) \, L_{n}(dy) -\frac{M}{p(n)}.
\end{equation}

Now let $A$ be a Borel set in $\mathcal M_1(\Sigma)$ and write
$$
\tilde{A} := \biggl\{ x  \in \Sigma^{p(n)} : L_n(x) \in A
\biggr\}.
$$
Using \eqref{Fn} and H\"older's inequality we obtain
\begin{eqnarray} \label{qv17}
P_n \bigl( L_{n} \in A \bigr) & = & \int_{\tilde{A}} \exp \biggl(
-\frac{2 n^2}{p(n)^2} \sum_{1 \leq i < j \leq p(n)}
F_n(x_i, x_j) \biggr) \nonumber \\
& & \exp \biggl(\frac{n}{p(n)} \sum_{i=1}^{p(n)} \log w_{n}(x_i)
\biggr)  {\rm m}_{\rr^{p(n)}}(dx)
\\
& \leq & \biggl(  \int  \exp \biggl(\frac{2 n}{p(n)} \log w_{n}(t)
\biggr) \, {\rm m}_{\rr}(dt) \biggr)^{\frac{p(n)}{{\tt 2}}} \nonumber \\
& &\left( \int_{\tilde{A}} \exp \left( - \frac{2 n^2}{p(n)^2}
\sum_{i \neq j} F_n^M(x_i, x_j) \right) {\rm m}_{\rr^{p(n)}}(dx)
\right)^{1/2}
 \nonumber\\
&=& {\rm (I)} \times {\rm (II)}. \nonumber
\end{eqnarray}
Note that we have used that $2 \sum_{i < j} F_n(x_i, x_j)  =
\sum_{i \neq j} F_n(x_i, x_j) $ by symmetry of $F_n$ in its
arguments.
 Now,
$\lim_{n \to \infty} \frac{1}{n^2} \log {\rm (I)} = 0$ by (a2). On
the other hand, for any $M >0$ \eqref{vor} yields
$${\rm (II)} \le \left\{ \int_{\tilde{A}} \exp \left( - 2 n^2 \left( L_{n}(x)^{\otimes 2}(F_n^M) - \frac{M}{p(n)} \right) \right) {\rm m}_{\rr^{p(n)}}(dx) \right\}^{1/2}$$
$$ \le \left\{ \exp \left( - 2 n^2 \left( \inf_{\mu \in A} \mu^{\otimes 2}(F_n^M) - \frac{M}{p(n)} \right) \right) \right\}^{1/2}$$
$$= \exp\left(- n^2  \inf_{\mu \in A} \mu^{\otimes 2}(F_n^M)\right) \exp\left( \frac{M n^2}{p(n)}\right).
$$

Using Lemma \ref{qv16}, we obtain
$$
\lim_{n \to \infty} \biggl( \inf_{\mu \in A} \mu^{\otimes
2}(F_n^M) \biggr) = \inf_{\mu \in A} \mu^{\otimes 2}(F^M).
$$
We have thus shown that for any Borel set $A \subset \mathcal
M_1(\Sigma)$ one has
\begin{equation} \label{qv18}
\limsup_{n \to \infty} \frac{1}{n^2} \log P_n \bigl( L_{n} \in A
\bigr) \leq - \inf_{\mu \in A} \int \int F^M(x,y) \, \mu(dx)
\mu(dy).
\end{equation}

Setting
$$ H(\mu) := \int F d\mu^{\otimes 2},\ H^M(\mu) := \int F^M d\mu^{\otimes
2},$$ one obtains well defined maps on $\calm_1(\Sigma)$ (see
Lemma \ref{qv16}). We will show that $H$ is a good rate function
that governs the LDP for $(P_n)$. To this end, observe that, $F^M$
being bounded and continuous, $H^M$ is weakly continuous on
$\mathcal M_1(\Sigma)$ for each $M
>0$. By monotone convergence, we have $\lim_{M \to \infty} H^M =
H$ pointwise on $\calm_1(\Sigma)$. As a limit of an increasing
sequence of continuous functions, $H$ is lower semi-continuous,
i.e. the level sets $\{H \leq L \}$ are closed. We claim that they
are compact. Indeed, let $m_F := |\inf F|$ and $a > 0$. Then for
any $\mu \in {\mathcal M}_1(\Sigma)$ one has
\begin{eqnarray*}
\biggl( \inf_{x,y \in [-a,a]^c} (F +m_F)(x,y) \biggr) \mu \bigl(
      [-a,a]^c \bigr)^2 & \leq & \int \int (F +m_F)(x,y) \, \mu(dx) \,  \mu(dy) \\
& \leq & H(\mu) +m_F,
\end{eqnarray*}
hence $\{H \leq L\} \subset K_L$, $L \in ]0, \infty[$, with
$$
K_L := \bigcap_{a>0} \biggl\{ \mu \in \mathcal M_1(\Sigma): \mu
\bigl( [-a,a]^c \bigr) \leq \biggl( \frac{L+m_F}{\inf_{x,y
\in[-a,a]^c} (F+m_F)(x,y)} \biggr)^{1/2} \biggr\}.
$$
Since $\lim_{a \to \infty} \inf_{x,y \in[-a,a]^c } (F+m_F)(x,y) =
\infty$, $K_L$ is weakly relatively compact by Prohorov's theorem. Hence the
rate function $H$ is good. Furthermore, this argument easily
yields the exponential tightness of $(P_n \circ L_n^{-1})_n$. Fix $M
> 0$, and define $K_L\ (L > 0)$ as above, using $F^M$ in the place
of $F$. For every $\mu \in K_L^c $ there exists $a = a_{\mu} > 0$
such that
$$ \left( \inf_{x, y \in [-a, a]^c} F^M(x, y) +
m_{F^M}\right)\ \mu([-a, a]^c)^2 > L + m_{F^M},$$ hence
\begin{equation}
\label{qv102} \inf_{x, y \in [-a, a]^c} F^M(x, y)\ \mu( [-a,
a]^c)^2 > L.\end{equation} Then \eqref{qv18} implies
\begin{eqnarray*}
 && \limsup_{n \to \infty} \frac{1}{n^2} \log
P_n(L_n \in K_L^c) \le - \inf_{\mu \in K_L^c} \int F^M
d\mu^{\otimes 2}\\
&=& - \inf_{\mu \in K_L^c} \left( \int_{\left([-a_{\mu},
a_{\mu}]^c\right)^2} F^M d\mu^{\otimes 2} + \int_{\rr^2 \setminus
\left([-a_{\mu}, a_{\mu}]^c\right)^2} F^M d\mu^{\otimes
2}\right)\\
&\le& - \inf_{\mu \in K_L^c} \int_{\left([-a_{\mu},
a_{\mu}]^c\right)^2} \inf_{x, y \in [-a_{\mu}, a_{\mu}]^c}F^M(x,
y)\ d \mu^{\otimes 2} - \inf_{\mu \in K_L^c} \int_{\rr^2 \setminus
\left([-a_{\mu}, a_{\mu}]^c\right)^2} \inf F^M d\mu^{\otimes 2}\\
&\le& - L + m_{F^M}.
\end{eqnarray*}
Since $\inf F^M > - \infty$, we have shown that
$$ \limsup_{L \to \infty} \limsup_{n \to \infty} \frac{1}{n^2} \log P_n(L_n \in K_L^c) =
- \infty,$$ hence that $(P_n \circ L_n^{-1})_n$ is exponentially
tight. \\

Now let $B(\mu, \delta)$ denote the ball centered at $\mu \in
\mathcal M_1(\Sigma)$ with radius $\delta$ for a distance
compatible with the weak topology. Since $\mu \mapsto H^M(\mu)$ is
continuous, from \eqref{qv18} we obtain for any $\mu \in \mathcal
M_1(\Sigma)$
$$
\inf_{\delta \to 0} \limsup_{n \to \infty} \frac{1}{n^2} \log P_n
\bigl( L_{n} \in B(\mu, \delta) \bigr) \leq - \int \int F^M(x,y) \,
\mu(dx)  \mu(dy).
$$
Finally, letting $M$ go to infinity, we obtain the following upper
bound
$$
\inf_{\delta \to 0} \limsup_{n \to \infty} \frac{1}{n^2} \log P_n
\bigl( L_{n} \in B(\mu, \delta) \bigr) \leq - H(\mu).
$$
\\

 Turning to the  {\it lower bound} for $(P_n \circ L_n^{-1})$, we show
that for any $\mu \in \calm_1(\Sigma)$
\begin{equation}
\label{qv19} \inf_{\delta > 0}   \liminf_{n \to \infty}
\frac{1}{n^2} \log P_n( L_n \in B(\mu, \delta))  \ge
\frac{\beta}{2} \kappa^2 \int \int \log |x^{\gamma} - y^{\gamma}|
\mu(dx)\mu(dy) + \kappa \int \log w d\mu.
\end{equation}
\begin{claim}
\label{qv20} We may assume, without loss of generality, that
\begin{itemize}
\item[(i)] $\mu$ has no atoms \item[(ii)] $\cals := {\rm
supp}(\mu)$ is a compact subset of $\Sigma$ such
that $\cals \cap (\caln(w) \cup \{ 0 \}) = \emptyset$.
\end{itemize}
\end{claim}
\begin{proof}  With the notations of the proof of the
upper bound, we know that
$$\frac{\beta}{2} \kappa^2 \int \int \log |x^{\gamma} - y^{\gamma}| \mu(dx)\mu(dy)
+ \kappa \int \log w(x) \mu(dx) = - \int F d\mu^{\otimes 2}$$ with
$F$ bounded from below. If $\mu$ has an atom, then $\int F
d\mu^{\otimes 2}$ is infinite, and there is nothing to prove. Set
$$ A_k := [-k, k] \cap \Sigma \cap \left( \bigcup_{x \in \caln \cup \{ 0 \}} ]x - \frac{1}{k}, x + \frac{1}{k}[ \right)^c$$
and $\mu_k := \frac{1}{\mu(A_k)} \mu|_{A_k}$. Then
$$ \int F d\mu^{\otimes 2} = \lim_{k \to \infty} \int F d\mu_k^{\otimes 2}.$$ Hence it suffices to prove (\ref{qv19}) for $\mu_k$ in the place of $\mu$. Consequently,
we may assume the support of $\mu$ to be contained in a finite
union of compact intervals not meeting $\caln(w) \cup \{ 0 \}$.
This implies (ii).
\end{proof}

For $j = 1, \ldots, p(n)$ let $\xi_j = \xi_j^{(n)}$ be the $\frac{
p(n) + 1 - j}{p(n)}$  quantile of $\mu$. Let $\gamma \in \nn$ be
as in \eqref{qv7} above. Write $\xi^{\gamma} =
(\xi_{p(n)}^{\gamma}, \ldots, \xi_1^{\gamma}).$ Set $\xi_{p(n) +
1} := \inf \cals$ and $\xi_0 := \xi_1 + 1$.
 Then, by Claim \ref{qv20},
$$ -  \infty < \xi_{p(n) + 1} < \xi_{p(n)} < \ldots < \xi_0 < \infty.$$

For $\delta > 0$, $t \in \rr^{p(n)}$ write
\begin{itemize}
\item $\pi_n(t) := \{ i = 1, \ldots, p(n):\ t_i \ge 0\},\ \nu_n(t)
:= \{ 1, \ldots, p(n)\} \setminus \pi_n(t)$,
\item $I_n(\delta) :=  \{ i = 1, \ldots, p(n):\ |\xi_i^{(n)} - \xi_{i+1}^{(n)}| \le \delta \}$,
\item
$\cali_j^{(n)}(t, \delta) := [t_j - \delta, t_j + \delta] \cap \Sigma,\ j = 1,
\ldots, p(n)$, \item $$ \calj_j^{(n)}(t, \delta) := \left\{
\begin{array}{ll} ~[t_j,  t_j + \delta]& {\rm  for}\ j \in \pi_n(t),\\
~[t_j - \delta,  t_j ] & {\rm  for}\ j \in \nu_n(t),
\end{array}\right. $$
\item ${\mathbb I}_n(t, \delta) := \prod_{j=1}^{p(n)}
\cali_j^{(n)}(t, \delta),$ ${\mathbb I}_n(t, \delta)^{\gamma} :=
\prod_{j=1}^{p(n)} \cali_j^{(n)}(t, \delta)^{\gamma},$ \item
${\mathbb J}_n(t, \delta) := \prod_{j=1}^{p(n)} \calj_j^{(n)}(t,
\delta).$
\end{itemize}
Here, for $M \subseteq \rr$, we write $M^{\gamma} := \{
m^{\gamma}:\ m \in M\}$.\\

 Fix $\delta > 0$, and write $G :=
B(\mu, 2\delta)$.
It follows from \cite[Lemma 3.3]{BAG} that for $n$ large enough
 one has
\begin{equation}
\label{delta.und.n}
      {\mathbb I}_n(\xi^{(n)}, \delta) \subset \{ x
\in \Sigma^{p(n)}:\ L_n(x) \in G \}.
\end{equation}

Set $\phi_j^{(n)} := \phi_j^{(n, \delta)} := \inf\{ w_n(x):\ x \in
[\xi_{j} - \delta , \xi_{j} + \delta] \cup [\xi_{j+1}, \xi_{j-1}]
\}, \ j = 1, \ldots, p(n)$,
 and write $\psi_n$ for the step function which equals $\phi_j^{(n)}$ on
$]\xi_{j+1}, \xi_j]$ and is zero elsewhere. \\

For $n$ as in \eqref{delta.und.n} we have
\begin{eqnarray*}
P_n ( L_n \in G ) &\ge& Z_n Q_n ({\mathbb I}_n(\xi, \delta))
= \int_{{\mathbb I}_n(\xi, \delta)} \prod_{i < j} |x_i^{\gamma} - x_j^{\gamma}|^{\beta}\ \prod_i w_n(x_i)^n {\rm m}_{p(n)}(dx)\\
&=& \frac{1}{\gamma^{p(n)}} \int_{{\mathbb I}_n(\xi,
\delta)^{\gamma}} \prod_{i < j} |x_i - x_j|^{\beta}\ \prod_i
w_n(x_i^{1 / \gamma})^n |x_i|^{- \frac{\gamma - 1}{\gamma}} {\rm
m}_{p(n)}(dx).
\end{eqnarray*}
Observing that $ {\mathbb I}_n(\xi, \delta)^{\gamma} \supset
{\mathbb J}_n(\xi^{\gamma}, \delta^{\gamma})$ and
$$ \inf \left\{ \prod_i w_n(x_i^{1 / \gamma})^n:\ x \in {\mathbb I}_{n}(\xi, \delta)^{\gamma} \right\}
=  \inf \left\{ \prod_i w_n(x_i)^n:\ x \in {\mathbb I}_{n}(\xi,
\delta) \right\} \ge \prod_i (\phi_i^{(n)})^n$$ we obtain that
$$P_n ( L_n \in G ) \ge
\frac{1}{\gamma^{p(n)}} \prod_i \left(\left( \phi_i^{(n)}\right)^n
(|\xi_i| + \delta)^{1 - \gamma}\right) \int_{ {\mathbb
J}_n(\xi^{\gamma}, \delta^{ \gamma})} \prod_{i < j}|x_i -
x_j|^{\beta} {\rm m}_{p(n)}(dx).$$ Now, we can bound this last
integral from below by
\begin{eqnarray}
&& \int_{ \left( [0, \delta^{\gamma}]^{\# \pi_n(\xi^{\gamma})}
\times [- \delta^{\gamma}, 0]^{\# \nu_n(\xi^{\gamma})} \right)
\cap \rr^{p(n), +}}
\prod_{i < j} | (\xi_i^{\gamma} + x_i) - (\xi_j^{\gamma} + x_j)|^{\beta}\ {\rm m}_{p(n)}(dx)\nonumber\\
&=&  \int_{ \left( [0, \delta^{\gamma}]^{\# \pi_n(\xi^{\gamma})}
\times [- \delta^{\gamma}, 0]^{\# \nu_n(\xi^{\gamma})} \right)
 \cap \rr^{p(n), +}}
\prod_{i < j} | (\xi_i^{\gamma} - \xi_j^{\gamma}) + (x_i - x_j)|^{\beta}\ {\rm m}_{p(n)}(dx)\nonumber\\
&\ge& \prod_{i < j-1} |\xi_i^{\gamma} - \xi_j^{\gamma}|^{\beta}\ \prod_{i=1}^{p(n) - 1} |\xi_{i+1}^{\gamma} - \xi_i^{\gamma}|^{\frac{\beta}{2}}\nonumber\\
\label{qv21} &\times&
 \int_{ \left( [0, \delta^{\gamma}]^{\# \pi_n(\xi^{\gamma})} \times [- \delta^{\gamma}, 0]^{\# \nu_n(\xi^{\gamma})} \right)
 \cap \rr^{p(n), +}}
\prod_{i=1}^{p(n) - 1} | x_{i} - x_{i+1}|^{\frac{\beta}{2}}\ {\rm
m}_{p(n)}(dx).
\end{eqnarray}
By the change of variables $u_{p(n)} = x_{p(n)},\ u_{i-1} = x_{i}
- x_{i-1}\ (i = p(n), \ldots, 2)$ one can bound (\ref{qv21}) from
below by
$$\int_{\left[0, \frac{\delta^{\gamma}}{ p(n)}\right]^{p(n)}} \prod_{i = 2}^{p(n)} u_i^{\frac{\beta}{2}}\ {\rm m}_{p(n)}(du)
=  \left( \frac{2}{\beta + 2} \right)^{p(n) - 1}\ \left(
\frac{\delta^{\gamma}}{ p(n)}\right)^{\frac{\beta + 2}{ 2} (p(n) -
1) + 1}.
$$
So far, it has been shown that
\begin{eqnarray*}
&&P_n ( L_n \in G)\\
&\ge&  \prod_{i < j-1} |\xi_i^{\gamma} - \xi_j^{\gamma}|^{\beta}\
\prod_{i=1}^{p(n) - 1} |\xi_{i}^{\gamma}
- \xi_{i+1}^{\gamma}|^{\frac{\beta}{2}} \\
&\times&
 \frac{1}{\gamma^{p(n)}} \prod_i \left( \left( \phi_i^{(n)}\right)^n (|\xi_i| + \delta)^{1 - \gamma}\right)
 \left( \frac{2}{\beta + 2} \right)^{p(n) - 1}\ \left( \frac{\delta^{\gamma}}{ p(n)}\right)^{\frac{\beta + 2}{ 2} (p(n) - 1) + 1}.
\end{eqnarray*}

So we obtain
\begin{eqnarray}
&&\frac{1}{n^2} \log P_n(L_n \in G)\nonumber\\
\label{term1.1}
&\ge& \left( \frac{p(n)}{n}\right)^2 \frac{\beta}{p(n)^2} \sum_{i < j} \log | \xi_i^{\gamma} - \xi_{j+1}^{\gamma}|\\
\label{term1.2}
&+& \left( \frac{p(n)}{n}\right)^2 \frac{\beta}{2 p(n)^2} \sum_{i = 1}^{p(n) -1} \log | \xi_{i+1}^{\gamma} - \xi_i^{\gamma}| \\
\label{term2}
&+& \frac{1}{n} \sum_i \log \phi_i^{(n)}\\
\label{term3}
&+& \frac{1}{n^2} \sum_i (1 - \gamma) \log (|\xi_i| + \delta) \\
\label{term4}
&+& \frac{1}{n^2}\left( - p(n)\log \gamma + (p(n) - 1) \log \frac{2}{\beta + 2}\right)\\
\label{term5} &+& \frac{1}{n^2} \left(\frac{(\beta + 2)(p(n) -
1)}{2} + 1\right) (\log (\delta^{\gamma}) - \log(p(n))).
\end{eqnarray}

Now, for (\ref{term1.1}) and (\ref{term1.2})   observe that
\begin{eqnarray*}
&& \int_{x < y} \log(y^{\gamma} - x^{\gamma}) \mu(dx) \mu(dy) \\
&=& \sum_{i < j} \int_{(x, y) \in [\xi_{j+1}, \xi_j] \times [\xi_{i+1}, \xi_i]} \log(y^{\gamma} - x^{\gamma}) \mu(dx) \mu(dy) \\
&+& \frac{1}{2} \sum_i \int_{(x, y) \in [\xi_{i+1}, \xi_i]^{\times 2}} \log|y^{\gamma} - x^{\gamma}| \mu(dx) \mu(dy) \\
&\le& \frac{1}{p(n)^2} \sum_{i<j} \log(\xi_i^{\gamma} -
\xi_{j+1}^{\gamma}) + \frac{1}{2 p(n)^2} \sum_i \log
(\xi_i^{\gamma} - \xi_{i+1}^{\gamma}).
\end{eqnarray*}
(\ref{term3}), (\ref{term4}), (\ref{term5}) are easily seen to
converge to zero.\\
As to \eqref{term2}, observe that $ \frac{1}{n} \sum_{j=1}^{p(n)}
\log \phi_j^{(n)} = \frac{p(n)}{n} \int \log \psi_n d\mu$. Denote
by $l$ a Lipschitz constant of $\log w$ on $\cals$. For $\eta > 0$
write $l_n^{\eta} := \max\{|\log w_n(x) - \log w_n(y)|:\ |x-y| \le
\eta\}$ and define $l^{\eta}$ analogously. By (a1), then,
$l_n^{\eta} \to l^{\eta}$. Note that $l^{\eta} \le l \eta$. Write
$M(n, \delta) := \bigcup_{j \in I_n(\delta)} [\xi_{j+1}^{(n)},
\xi_j^{(n)}]$ and $C := \max\{ \log w(x):\  x \in \cals \} -
\min\{ \log w(x):\ x \in \cals \}.$

Let $\epsilon > 0$. Since for all $n$ and $j \ge 1$ one has
$\mu\left([\xi_{j+1}^{(n)}, \xi_j^{(n)}]\right) = 1 / p(n)$, and
since $\delta $ is fixed, one has $\mu( M(n, \delta)^c) \le
\epsilon$ for large $n$. Now let $n$ be large enough such that one
also has $|l_n^{2\delta} - l^{2\delta}| \le \epsilon$ and $\| \log
w_n - \log w\| \le \epsilon$.  Then
\begin{eqnarray*}
 \int |\log \psi_n - \log w| d\mu &\le&  \int |\log \psi_n - \log w_n| d\mu + \epsilon\\
&\le& \sum_{j \in I_n(\delta)} \int_{[\xi_{j+1}, \xi_j]} |\log \psi_n - \log w_n | d\mu + (C+ 2\epsilon +1) \epsilon \\
&\le& p(n) \frac{1}{p(n)}\ 2l\delta +  (C+ 2\epsilon +3) \epsilon.
\end{eqnarray*}
So we have $$ \limsup_{n \to \infty} \left| \frac{1}{n}
\sum_{j=1}^{p(n)} \log \phi_j^{(n)} - \kappa \int \log w\
d\mu\right| = O(\delta).$$
This implies \eqref{qv19}.\\

Summing up, for any $\mu \in \calm_1(\Sigma)$ we have obtained
$$ \inf_{\delta > 0} \limsup_{n \to \infty} \frac{1}{n^2}
\log P_n(L_n \in B(\mu, \delta)) \le - H(\mu)$$ and
$$ \inf_{\delta > 0} \liminf_{n \to \infty}
\frac{1}{n^2} \log P_n(L_n \in B(\mu, \delta)) \ge - H(\mu).$$
Using the exponential tightness of $(P_n \circ L_n^{-1})_n$, we can apply
\cite[Thm.\ 4.1.11]{DZ} to obtain an LDP for $(P_n \circ L_n^{-1})_n$ with rate
$H$ and speed $n^2$. Setting $A = G = \calm_1(\Sigma)$ in the
lower and upper bound, we obtain
$$ \lim_{n \to \infty} \frac{1}{n^2} \log Z_n = - \inf_{\mu \in
\calm_1(\Sigma)} \int F d\mu^{\otimes 2}.$$ By Lemma \ref{qv16}
(ii), the right-hand side is $< + \infty$. Now,
$$ \frac{1}{n^2} \log Q_n(L_n \in A) = \frac{1}{n^2} \left( \log
P_n(L_n \in A) - \log Z_n\right)$$ for any Borel set $A$ in
$\calm_1(\Sigma)$. Hence Theorem \ref{qv10} is proven.

\section{Application to Gaussian Hamiltonian ensembles}
\label{qv43}

\subsection{Gaussian Hamiltonian ensembles}

 Let $\mm^{(n)}_{\calc} \subset \cc^{d(n) \times
d(n)}$ be the space of good Hamiltonians of symmetry  class
$\calc$. We wish to define ${\rm Ad}_G$- resp.\ ${\rm
Ad}_K$-invariant probability measures on $\mm^{(n)}_{\calc}$ with
the additional property that all matrix entries should be Gaussian
and as many entries as possible should be independent. Plainly,
what we are looking for is an analog for class $\calc$ of the GOE,
which was reviewed in
the introduction and corresponds to $\calc =$ AI. \\

Let us look at an example: For $\calc =$ CI, represent an
element of $\mm_{{\rm CI}}^{(n)} $ as
$$ X = \left( \begin{array}{rr} A & B \\ B & -A \end{array}
\right),$$ where $A$ and $B$ are real symmetric $n \times
n$ matrices. The upper diagonal entries of $A$ and $B$ can be chosen
as independent Gaussians. Then the
distribution of $X$ has a Lebesgue density of the form
$$ {\rm const} \times \prod_{k \le l} \exp\left(-  \frac{ a_{kl}^2}{2
\sigma_{a,k,l}^2}\right) \prod_{k \le l} \exp\left(-  \frac{
b_{kl}^2}{2 \sigma_{b,k,l}^2}\right).$$ Now, a suitable choice of
$\sigma_{a, k, l}, \sigma_{b, k, l}$ guarantees the invariance of
the distribution of $X$. In fact, $X$ being symmetric,

\begin{eqnarray*}
&& {\rm Tr}(X^2) = {\rm Tr}(X X\tra) = 2\ {\rm Tr}(A A\tra + B
B\tra) = 2 \sum_{i, j} a_{ij}^2 +
  b_{ij}^2\\
&=& 4 \sum_{i < j} a_{ij}^2 + 2 \sum_i a_{ii}^2  + 4 \sum_{i < j} b_{ij}^2 + 2 \sum_i b_{ii}^2,
\end{eqnarray*}
hence
$$
- \frac{1}{8 \sigma^2} {\rm Tr}(X^2) =  \sum_{i < j} - \frac{a_{ij}^2}{2
  \sigma^2} +
\sum_i - \frac{a_{ii}^2}{2(2\sigma^2)} +  \sum_{i < j} - \frac{b_{ij}^2}{2 \sigma^2}
+ \sum_i - \frac{b_{ii}^2}{2(2\sigma^2)}
$$

This means that in order to obtain an invariant distribution, we
choose $\sigma^2 > 0$ and set $ \sigma_{a, k, l}^2 = \sigma_{b, k,
l} = \sigma^2\ (k < l),\ \sigma_{a, k, k} = \sigma_{b, k, k}  = 2
\sigma^2.$ We write ${\rm GE}_{{\rm CI}}^{(n)}(\sigma^2)$ for the
probability distribution on $\mm_{{\rm CI}}^{(n)}$ obtained in
this way.

Note that the eigenvalues of $X$ come in pairs $\pm \lambda$. So, if
$\lambda_1, \ldots, \lambda_n$ are the positive eigenvalues of $X$ (under
${\rm GE}_{{\rm CI}}^{(n)}(\sigma^2)$, $X$ has $n$ distinct positive
eigenvalues a.s.), then we have
$$ - \frac{1}{8 \sigma^2} {\rm Tr}(X^2) = - \frac{1}{8 \sigma^2}\ 2\
\sum_{j=1}^n \lambda_j^2 =
- \frac{1}{4 \sigma^2}\
\sum_{j=1}^n \lambda_j^2.$$

For comparison with the other symmetry classes, we write $\phi_{{\rm CI}} = 8$
and $\psi_{{\rm CI}} = 4$. \\

For a general symmetry class $\calc$, $X \in \mm_{\calc}^{(n)}$,
we may proceed along the same lines. If the nonzero real and
imaginary parts of the strictly upper triangular entries of (skew)
symmetric or hermitian blocks have variance $\sigma^2$, then the
nonzero real and imaginary parts of the diagonal entries of these
blocks have variance $2\sigma^2$. This procedure determines an
integer $\phi_{\calc}$ such that the Lebesgue density of $X$ has
the form
$$ {\rm const.} \times \exp\left( - \frac{1}{\phi_{\calc} \sigma^2} {\rm Tr}(X^2)\right).$$
If $\lambda_1, \ldots, \lambda_{p(n)}$ are the positive
eigenvalues without multiplicity of $X$ (or all eigenvalues
without multiplicity for $\calc =$ A, AI, AII), then there is an
integer $\psi_{\calc}$ defined by
$$ - \frac{1}{\phi_{\calc} \sigma^2} {\rm Tr}(X^2) = -  \frac{1}{\psi_{\calc}
  \sigma^2} \sum_{j = 1}^{p(n)} \lambda_j^2.$$
Proceeding analogously to the example above and keeping in mind
that $X$ has $2$-dimensional eigenspaces for
  $\calc =$ AII, CII and DIII, one obtains the
  following data:

\begin{center}
\begin{tabular}{|l|l|l|l|l|l|l|l|l|l|l|}
\hline
$\calc$ & A & AI & AII & AIII & B/D & BDI & DIII & C & CI & CII\\
\hline $\phi_{\calc} $ & 4 & 4 & 8 & 4 & 4 & 4 & 8 & 8 & 8 & 8 \\
\hline $ \psi_{\calc}$ & 4 & 4 & 4 & 2 & 2 & 2 & 2 & 4 & 4 & 2\\
\hline
\end{tabular}
\end{center}

In the following subsections we apply our main theorem to the ensembles
${\rm GE}_{\calc}^{(n)}(\frac{\sigma^2}{n})$.

\subsection{LDP for Wigner-Dyson ensembles A, AI, AII}

From Corollary \ref{qv6} we see that one can subsume the joint
eigenvalue density induced by ${\rm GE}_{\calc}^{(n)} \left(
\frac{\sigma^2}{n}\right)$, $\calc =$ A, AI, AII, under the
general form of \eqref{qv7} by choosing $p(n) = n$, hence $\kappa
= 1$, $\gamma = 1$, $\beta = 1, 2, 4$ according to $\calc =$ AI,
A, AII, and
\begin{equation}
\label{qv26} w_n(x_j) = \exp\left( - \frac{1}{4 \sigma^2}\
x_j^2\right),
\end{equation}
independent of $n$. Then, under ${\rm GE}_{\calc}^{(n)} \left(
\frac{\sigma^2}{n}\right)$, $(L_n)_n$ satisfies an LDP with good
rate function
\begin{equation}
\label{qv27} I_{\calc}(\mu)  = I_{\calc, \sigma^2}(\mu) =
\frac{\beta}{2} \int \int \log |x-y|^{-1} \mu(dx) \mu(dy) +
\frac{1}{4 \sigma^2}\ \int x^2 \mu(dx) - {\rm const}.
\end{equation}
To apply the theory of logarithmic potentials with external
fields, as presented in \cite{ST}, we exploit the fact that the
unique minimizer $\mu^*$ of $I_{\calc}$ is also the unique
minimizer of
$$ J_{\calc}(\mu) =
 \int \int \log |x-y|^{-1} \mu(dx) \mu(dy) +
 2 \int \frac{1}{4 \sigma^2 \beta} x^2 \mu(dx).
$$
One reads off from p.\ 284 (cf.\ p.\ 26 for notations) of \cite{ST} that $\mu^*$
has Lebesgue density
$$ 1_{[- 2 \sqrt{\sigma^2 \beta},\ 2 \sqrt{\sigma^2 \beta}]}(x)\ \frac{1}{2 \pi \sigma^2 \beta}
\sqrt{4 \sigma^2 \beta - x^2},$$
the density of Wigner's {\it semicircle distribution} of radius
$2 \sqrt{\sigma^2 \beta}$.

\subsection{LDP for chiral ensembles BDI, AIII, CII}
\label{qv66} We will freely use the notation introduced above and
in Section \ref{qv2}. We have $\psi_{\calc} = 2$ and $p(n) = s(n)
\wedge t(n)$. For simplicity, we will assume that $s(n) \le t(n)$
for all $n \in \nn$, hence $\alpha(n) = \beta(t(n) - s(n)) + \beta
- 1$. Then we can subsume the joint eigenvalue density induced by
${\rm GE}_{\calc}^{(n)} \left( \frac{\sigma^2}{n}\right)$ under
the general form of \eqref{qv7} by setting $\gamma = 2$ and
$$ w_n(x) = x^{\frac{\beta(t(n) - s(n)) + \beta - 1}{n}} e^{- \frac{x^2}{2 \sigma^2}},$$ hence
$$ w(x) = x^{\beta(1 - 2 \kappa)} e^{- \frac{x^2}{2 \sigma^2}}.$$
Then, by Theorem \ref{qv10}, the rate function of the LDP for $Q_n
\circ L_n^{-1}$ is
\begin{equation}
\label{qv55}
 I_{\calc}(\mu) = \beta \frac{\kappa^2}{2} \int \int \log
\frac{1}{|x^2 -y^2|} \mu(dx) \mu(dy) - \kappa   \int  \log
\left(x^{ \beta(1 - 2\kappa)}   e^{- \frac{x^2}{2
\sigma^2}}\right) \mu(dx) - c,
\end{equation}
for $\mu \in \calm_1([0, \infty[).$ For $r > 0, x \ge 0$ write
$T_r(x) := x^r$. Now, $\mu^*$ is the unique minimizer of
$I_{\calc}$ if, and only if, $T_2(\mu^*)$ is the unique minimizer
of
$$ \beta \frac{\kappa^2}{2}  \int \int  \log
\frac{1}{|x -y |} \nu(dx) \nu(dy) -  \kappa   \int  \log
\left(x^{\frac{ \beta(1 - 2\kappa)}{2}}   e^{- \frac{x}{2
\sigma^2}}\right) \nu(dx) - c,$$ hence of
\begin{equation}
\label{qv56} J_{\calc}(\nu) := \int \int \log \frac{1}{|x -y |}
\nu(dx) \nu(dy) -  2  \int  \log \left(x^{\frac{ 1 - 2\kappa}{2
 \kappa}}   e^{- \frac{x}{2 \beta \kappa \sigma^2}}\right)
\nu(dx).
\end{equation}
We can use the following facts from logarithmic potential theory
with Laguerre weights:
\begin{lem}
\label{qv57} For $s \ge 0,\ \lambda > 0$ the integral $$I_{s,
\lambda}(\nu) := \int \int  \log \frac{1}{|x -y |} \nu(dx) \nu(dy)
-  2 \int \log \left( x^s e^{- \lambda x}\right) \nu(dx),\quad \nu
\in \calm_1([0, \infty[)$$ has a unique minimizer $\nu^*_{s,
\lambda}$ with Lebesgue density
\begin{equation}
\label{qv59}
 1_{[a, b]}(x)\ \frac{\lambda}{\pi x} \sqrt{(x-a)(b-x)},
 \end{equation}
where
\begin{equation}
\label{qv58}
 a = a_{s, \lambda} = \frac{1}{\lambda} (s + 1 - \sqrt{ 2s +
1}),\  b = b_{s, \lambda} = \frac{1}{\lambda} (s + 1 + \sqrt{ 2s +
1}).
\end{equation}
Consequently, $T_{\frac12}(\nu^*_{s, \lambda})$ has Lebesgue
density
\begin{equation}
\label{qv65} 1_{[ \sqrt{a}, \sqrt{b}]}(x)\ \frac{2 \lambda}{\pi x}
 \sqrt{ (x^2 - a) (b - x^2)}.
 \end{equation}
\end{lem}
\begin{proof} \cite{ST}, IV (1.31), IV (5.18) \end{proof}
We have  $s = \frac{1}{2\kappa} - 1$ and $\lambda = \frac{1}{2
\sigma^2 \beta \kappa}.$ Note that $2\kappa \le 1$ by our
assumptions. This yields
\begin{equation}
\label{qv60} a = 2 \sigma^2 \beta \left( \frac{1}{2} -
\sqrt{\kappa (1 - \kappa)}\right),\ b = 2 \sigma^2 \beta \left(
\frac{1}{2} + \sqrt{\kappa (1 - \kappa)}\right),
\end{equation}
and from \eqref{qv65} we conclude that the minimizer $\mu^*$ of
$I_{\calc}$ has Lebesgue density
\begin{equation}
\label{qv61}
 1_{[ \sqrt{a}, \sqrt{b}]}(x)\ \frac{1}{ \sigma^2 \beta \kappa \pi x}
 \sqrt{ (x^2 - a) (b - x^2)}
\end{equation}
with $a, b$ as in \eqref{qv60}.
\begin{remark}
\label{qv101}
$T_2(\mu^*)$ is related to, but does not coincide with, the Mar\v{c}enko-Pastur distribution, which arises as limit of the empirical eigenvalue measure of Wishart (or more generally, sample covariance) matrices. Details on this relationship can be found in the companion paper \cite{KHC}. The LDP for Wishart matrices, first proven by Hiai and Petz in \cite{HPwish}, is a consequence of Theorem \ref{qv10} above.
\end{remark}

\subsection{LDP for BdG ensembles B, D, DIII, C, CI}

Although there are five BdG ensembles, there are only four natural
large $n$ limits for these ensembles. This is because B and D are
$i \frs \fro(n)$ for odd resp.\ even $n$, and, as $n$ tends to
infinity, these series should be thought of as interlaced. It is
consistent to do so, because we will see that the parameter
$\alpha$, where B and D differ (see table in Corollary \ref{qv6}),
does not affect the rate function.  By the same token, it is in
fact consistent to consider one rather than two large $n$ limits for DIII, because the even and odd cases of DIII only differ w.r.t.\ $\alpha$.\\

We subsume the joint eigenvalue density induced by ${\rm
GE}_{\calc}^{(n)} \left( \frac{\sigma^2}{n}\right)$, $\calc =$ B,
D, DIII, C, CI under the general form \eqref{qv7} by setting  $\gamma = 2$, $(\alpha, \beta)$
according to the table after Corollary \ref{qv6}, $p(n) = \lfloor \frac{n}{2} \rfloor$ for $\calc =$
DIII and $p(n) = n$ otherwise, hence  $\kappa = \frac{1}{2}$ or $\kappa = 1$, respectively. This yields
\begin{equation*}
\label{qv28} w_n(x) = x^{\alpha/n} \exp( - \frac{1}{\psi_{\calc} \sigma^2}\
x^2),
\end{equation*} hence
\begin{equation*}
\label{qv29} w(x) = \exp( -  \frac{1}{\psi_{\calc} \sigma^2}\ x^2 )
\end{equation*}
with $\psi_{{\rm B/D}} = 2$ and $\psi_{{\rm C}} = \psi_{{\rm CI}} = \psi_{{\rm DIII}} = 4$.

Then, under ${\rm GE}_{\calc}^{(n)} \left(
\frac{\sigma^2}{n}\right)$, $(L_n)_n$ satisfies an LDP with good
rate function
\begin{equation*}
\label{qv30} I_{\calc}(\mu) = \frac{\beta}{2} \kappa^2 \int \int \log
|x^2-y^2|^{-1} \mu(dx) \mu(dy) + \frac{\kappa}{\psi_{\calc} \sigma^2}\ \int x^2
\mu(dx) - {\rm const}.
\end{equation*}
To determine the unique minimizer $\mu^*$ of $I_{\calc}$, one can
proceed as in Subsection \ref{qv66} and apply Lemma \ref{qv57}
with $s = 0$ and $\lambda = \frac{1}{\psi_{\calc} \sigma^2 \beta \kappa}$.
By \eqref{qv65},  $\mu^*$ is a {\it quarter circle distribution},
given by the Lebesgue density
$$ 1_{[ 0, \sqrt{2 \psi_{\calc} \sigma^2 \beta \kappa}]}(x)\ \frac{2}{\psi_{\calc} \sigma^2 \beta \kappa \pi}\
\sqrt{ 2 \psi_{\calc} \sigma^2 \beta \kappa - x^2}.$$

\bibliographystyle{amsplain}
\bibliography{az}

\end{document}